\documentclass{gtmon_a}
\pdfoutput=1
\usepackage[T1,T5]{fontenc}
\usepackage{rotating}

\proceedingstitle{Proceedings of the School and Conference in Algebraic
Topology (The Vietnam National University, Hanoi, 9-20 August 2004)}
\conferencestart{9 August 2004}
\conferenceend{20 August 2004}
\conferencename{School and Conference in Algebraic Topology}
\conferencelocation{Vietnam National University, Hanoi, Vietnam}

\editor{John Hubbuck}
\givenname{John}
\surname{Hubbuck}

\editor{Nguy\~\ecircumflex{}n H V H\uhorn{}ng}
\givenname{H\uhorn{}ng}
\surname{Nguy\~\ecircumflex{}n}

\editor{Lionel Schwartz}
\givenname{Lionel}
\surname{Schwartz}

\title{The Lambda algebra and $\mathrm{Sq}^0$}

\author{J\,H Palmieri}
\givenname{J\,H}
\surname{Palmieri}
\address{Department of Mathematics\\
University of Washington\\\newline
Box 354350\\
Seattle WA 98195-4350\\
USA}
\email{palmieri@math.washington.edu}
\urladdr{}

\volumenumber{11}
\issuenumber{}
\publicationyear{2007}
\papernumber{10}
\startpage{201}
\endpage{216}

\doi{}
\MR{}
\Zbl{}

\arxivreference{}


\received{1 April 2005}
\revised{8 October 2005}
\accepted{15 November 2005}
\published{14 November 2007}
\publishedonline{14 November 2007}
\proposed{}
\seconded{}
\corresponding{}
\editor{}
\version{}

\makeatletter
\def\cnewtheorem#1[#2]#3{\newtheorem{#1}{#3}[section]
\expandafter\let\csname c@#1\endcsname\c@thm}

\input xypic
\let\xysavmatrix\xymatrix
\def\xymatrix{\disablesubscriptcorrection\xysavmatrix}
\AtBeginDocument{\let\wbar\wbar\let\tilde\wtilde\let\hat\what}

\theoremstyle{plain}
\newtheorem{thm}{Theorem}[section]
\cnewtheorem{prop}[thm]{Proposition}
\cnewtheorem{cor}[thm]{Corollary}
\cnewtheorem{lemma}[thm]{Lemma}

\theoremstyle{definition}
\cnewtheorem{definition}[thm]{Definition}
\cnewtheorem{conj}[thm]{Conjecture}
\cnewtheorem{exmp}[thm]{Example}

\theoremstyle{remark}
\cnewtheorem{rem}[thm]{Remark}
\cnewtheorem{rems}[thm]{Remarks}
\cnewtheorem{note}[thm]{Note}
\cnewtheorem{warn}[thm]{Warning}
\newtheorem*{question}{Question}

\makeatother  

\DeclareMathOperator{\Ext}{Ext}
\DeclareMathOperator{\SSq}{Sq}
\DeclareMathOperator{\Span}{Span}
\DeclareMathOperator{\im}{im}
\newcommand{\Sq}{\SSq^{0}}
\newcommand{\ztwo}{\mathbb{F}_{2}}
\newcommand{\ext}[2]{\Ext_{#1}^{#2}(\ztwo, \ztwo)}
\newcommand{\lie}{A_{\mathsf{Lie}}}
\newcommand{\suchthat}{:}
\newcommand{\quotient}{\Lambda'}
\newcommand{\cocomplete}{\theta^{-1} \Lambda}
\newcommand{\complete}{\widehat{A}}
\newcommand{\zhalf}{{\textstyle \Z[\frac{1}{2}]}}
\newcommand{\odd}{\text{odd}}
\newcommand{\even}{\text{even}}

\renewcommand{\theenumi}{\alph{enumi}}

\begin{document}

\begin{abstract} 
The action of $\mathrm{Sq}$ on the cohomology of the Steenrod algebra is
induced by an endomorphism $\theta$ of the Lambda algebra.  This paper
studies the behavior of $\theta$ in order to understand the action of
$\mathrm{Sq}$; the main result is that $\mathrm{Sq}$ is injective in
filtrations less than 4, and its kernel on the 4--line is computed.
\end{abstract}

\maketitle

\section{Introduction}

The Lambda algebra $\Lambda$ was constructed in the ``six author
paper'' \cite{6-author}, and numerous authors have studied it over the
past 40 years.  In \cite{wang.j*;hopf}, Wang studied an algebra
endomorphism $\theta$ of $\Lambda$, and this endomorphism is the main
focus of this paper.

$\Lambda$ is the bigraded differential $\ztwo$--algebra with generators
$\lambda_{n} \in \Lambda^{1,n+1}$ for $n \geq 0$, with relations
and differential given as follows:
\begin{gather*}
\sum_{i+j=n} \tbinom{i+j}{i} \lambda_{i-1+m} \lambda_{j-1+2m} \quad
\quad \textup{for} \ m \geq 1, n \geq 0, \\
d(\lambda_{n-1}) = \sum_{i+j=n} \tbinom{i+j}{j} \lambda_{i-1}
\lambda_{j-1}.
\end{gather*}
Let $A$ be the mod 2 Steenrod algebra.  The following result and its
unstable analogue are the main reasons for studying the Lambda algebra.

\begin{thm}[Bousfield et al.~\cite{6-author}]
\label{thm-6-author}
$H^{s,t}(\Lambda) \cong \ext{A}{s,t}$.  Indeed, one
can view $\Lambda$ as the $E_{1}$--term of the classical Adams spectral
sequence converging to the 2--component of the stable homotopy groups
of spheres.
\end{thm}

Given the free algebra $F = \textup{alg}_{\ztwo}(\lambda_{0},
\lambda_{1}, \dotsc)$ of which $\Lambda$ is a quotient, define a map
$\wtilde{\theta} \co F \rightarrow F$ by $\wtilde{\theta}
(\lambda_{n}) = \lambda_{2n+1}$.

\begin{prop}[Wang \cite{wang.j*;hopf}]
\label{prop-wang-lambda}
$\wtilde{\theta}$ induces a map $\theta \co \Lambda \rightarrow
\Lambda$ of differential graded algebras, and this map is one-to-one.
\end{prop}

\begin{cor}
The subalgebra of $\Lambda$ generated by $\{ \lambda_{2n+1} \suchthat
n \geq 0 \}$ is a sub-differential graded algebra, which is isomorphic
to $\Lambda$.
\end{cor}

Note that the isomorphism $\theta \co \Lambda \smash{\xrightarrow{\cong}}
\im (\theta)$ doubles internal degrees.  There is an endomorphism, called
$\Sq$, of $H^{*}\Lambda = \ext{A}{*}$ which also doubles internal
degrees.  One way to describe this map is as follows: the graded dual
of $A$ is a commutative Hopf algebra over $\ztwo$, and the Frobenius
map is a Hopf algebra map.  $\Sq$ is the induced map on $\ext{A}{*}$.
See May \cite[Proposition 11.10]{may;steenrod} for this result, and see
Palmieri
\cite{palmieri;quotient} for more discussion of $\Sq$ acting on
Ext.

\begin{prop}\label{prop-theta}
$\theta \co \Lambda \rightarrow \Lambda$ induces the map
\[
\Sq \co \ext{A}{*} \rightarrow \ext{A}{*}.
\]
\end{prop}

Bob Bruner told me how to prove this.

\begin{proof}
According to Priddy \cite{priddy;koszul}, $\lambda_{n}$ is
represented in the cobar construction by $[\xi_{1}^{n+1}]$.  According
to May \cite[Proposition 11.10]{may;steenrod}, $\Sq$ is induced by
the dual to the Frobenius, so $\Sq$ takes $[\xi_{1}^{n+1}]$ to
$[\xi_{1}^{2n+2}]$; thus it agrees with the map induced by $\theta$,
which sends $\lambda_{n}$ to $\lambda_{2n+1}$.
\end{proof}

The goal of this paper is to use the DGA endomorphism $\theta \co
\Lambda \rightarrow \Lambda$ to study the action of $\Sq$ on Ext over
the Steenrod algebra.  The main result is the following; the result
when $s \leq 3$ was first proved by Wang.

\begin{thm}\label{thm-main}
Consider $\Sq \co \ext{A}{s,*} \rightarrow
\ext{A}{s,*}$.  When $s \leq 3$, this map is
injective.  When $s=4$, the kernel is one-dimensional: 
\[
\ker (\Sq|\Ext_{A}^{4}) = \Span (h_{0}^{4}).
\]
\end{thm}

This theorem is proved in \fullref{sec-bss} using a Bockstein
spectral sequence argument.  The Bockstein spectral sequence is also
of interest in its own right.  We also invert $\theta$ to construct
the ``cocomplete Lambda algebra'' in \fullref{sec-cocomplete}.

This paper studies some of the same issues as
\cite{palmieri;quotient}, but uses a different approach.  As with that
paper, this work raises more questions than it answers; here are
some.

\begin{question}
\begin{enumerate}
\item Is $(\Sq)^{-1} \ext{A}{*}$ finite dimensional in each bidegree?
\item Even better, for fixed $s$, is there a uniform bound on 
\[
\dim (\Sq)^{-1} \ext{A}{s,t}?
\]
\item H\uhorn{}ng has asked the following: for fixed $s$, is there a bound
on the ``$\Sq$--nilpotence height'' of elements in
$\ext{A}{s,*}$?  That is, is there a number $N=N(s)$ so that, for
$z \in \ext{A}{s,*}$, either $(\Sq)^{N}(z)=0$ or
$(Sq^{0})^{i}(z)\neq 0$ for all $i$?  See Hu'ng
\cite{hung;cohomology-steenrod} for more details.
\end{enumerate}
\end{question}

The Bockstein spectral sequence of \fullref{sec-bss} may be
helpful in answering the last question, and perhaps the cocomplete
Lambda algebra (\fullref{sec-cocomplete}) could be helpful for the
first two.

\section{Recollections on the Lambda algebra}\label{sec-lambda}

In this section, we recall some facts about the Lambda algebra.
$\Lambda$ was first constructed in the ``six author paper''
\cite{6-author}.  Wang \cite{wang.j*;hopf} gave alternate forms for
the relations and the differential, described the endomorphism
$\theta$, and used Lambda algebra calculations to prove the Hopf
invariant one problem.  See also Ravenel's green book \cite[Section
3.3]{ravenel;green}, which summarizes the basic facts about $\Lambda$,
and also describes the Curtis algorithm; Priddy's Koszul resolutions
paper \cite{priddy;koszul} gives a purely algebraic construction of
$\Lambda$; and Richter has given elementary proofs of some of the
properties of $\Lambda$ in \cite{richter;lambda}.

$\Lambda$ is a bigraded differential $\ztwo$--algebra; it is defined to
be 
\[
\Lambda = \textup{alg}_{\ztwo}(\lambda_{0}, \lambda_{1}, \lambda_{2},
\dotsc) / (\text{relations}),
\]
graded by putting $\lambda_{n}$ in bidegree $(1,n+1)$.  There are two
forms for the relations and the differential, the ``symmetric'' form 
\begin{gather}\label{reln-symm}
\sum_{i+j=n} \tbinom{i+j}{i} \lambda_{i-1+m} \lambda_{j-1+2m} \quad
\quad \textup{for} \quad m \geq 1, n \geq 0, \\
d(\lambda_{n-1}) = \sum_{i+j=n} \tbinom{i+j}{j} \lambda_{i-1}
\lambda_{j-1},
\end{gather}
and the ``admissible'' form
\begin{gather}\label{reln-admissible}
\lambda_{i} \lambda_{2i+1+n} = \sum_{j \geq 0} \tbinom{n-j-1}{j}
\lambda_{i+n-j} \lambda_{2i+1+j} \quad \quad \textup{for} \quad i, n
\geq 0, \\
d(\lambda_{n}) = \sum_{j \geq 1} \tbinom{n-j}{j}
\lambda_{n-j}\lambda_{j-1}.
\end{gather}
Wang \cite{wang.j*;hopf} showed that these two sets of relations
generate the same ideal in the free algebra on the $\lambda_{n}$'s.
One can show that $d$ behaves formally like taking the commutator with
$\lambda_{-1}$; Bruner noted this in \cite{bruner;augmented}.

We refer to the first grading in $\Lambda$ as the
\emph{homological} grading, and the second as the \emph{internal}
grading.  Note that some authors define $\lambda_{n}$ to have degree
$(1,n)$ rather than $(1,n+1)$, in which case the second grading is the
\emph{stem} grading, the difference between the internal and
homological gradings.

\begin{definition}
A monomial $\lambda_{i_{1}} \dotsb \lambda_{i_{s}}$ is
\emph{admissible} if $2i_{r} \geq i_{r+1}$ for $1 \leq r \leq s-1$.
\end{definition}

\begin{prop}[Bousfield et al.~\cite{6-author}]
\label{prop-admissible}
The admissible monomials form a basis for $\Lambda$.
\end{prop}

Note that, while it is easy to show from the relations above that the
admissible monomials span $\Lambda$, it is not trivial to show that
they are linearly independent.  Similarly, it is not obvious that the
differential $d$ respects the relations.  See \cite{richter;lambda}
for elementary proofs.

As remarked in \fullref{thm-6-author}, the cohomology
$H^{*}(\Lambda,d)$ is $\Lambda$'s raison d'\^etre.

\section{A Bockstein spectral sequence}\label{sec-bss}

Let
\[
\quotient = \Lambda / \theta \Lambda,
\]
and consider the short exact sequence of chain complexes
\begin{equation}\label{eqn-ses}
0 \rightarrow \Lambda \xrightarrow{\theta} \Lambda \rightarrow
\quotient \rightarrow 0.
\end{equation}
After taking cohomology, this gives a Bockstein spectral sequence
with
\begin{gather}\label{eqn-bss}
E_{1}^{s,t,u} = 
\begin{cases}
H^{s+t,u}(\quotient), & \text{when} \ \ s \geq 0, \\
0 & \text{otherwise},
\end{cases} \\
d_{r} \co E_{r}^{s,t,u} \rightarrow E_{r}^{s+r,t-r+1,u/2^{r}},
\end{gather}
converging to $H^{s+t,u} (\Lambda) = \ext{A}{s+t,u}$.  The phrase
``Bockstein spectral sequence'' is perhaps ambiguous, so some details
may be helpful.  One way to construct it is to filter $\Lambda$ by
setting $F^{s}\Lambda = \im \theta^{s}$.  Because $\theta$ is
injective, $F^{s}\Lambda$ is isomorphic to $\Lambda$; more precisely,
$F^{s}\Lambda$ in bidegree $(i,2^{s}j)$ is isomorphic to $\Lambda$ in
bidegree $(i,j)$.  Up to a similar doubling of internal degrees,
$F^{s}\Lambda / F^{s+1}\Lambda$ is isomorphic to $\quotient$.  The
Bockstein spectral sequence is the one which arises from taking
cohomology of this filtered chain complex.  Thus for fixed $s \geq 0$
and $t$ we have $E_{1}^{s,t} \cong E_{1}^{s+i,t-i}$ for all $i$
with $i \geq 0$.  Given a class $x \in E_{1}^{s,0}$, we write
$\wbar{\theta}^{-i} x$ for the image of $x$ under this
isomorphism.  The differentials are ``$\wbar{\theta}$--periodic''
and determined by their effect on $E_{r}^{*,0}$. For any $x \in
E_{r}^{s,0}$, there is a differential
\[
d_{r} \co x \longmapsto \wbar{\theta}^{r} y,
\]
if and only if for any $i \geq 0$, there is a differential
\[
d_{r} \co \wbar{\theta}^{i} x \longmapsto \wbar{\theta}^{i+r} y.
\]
Furthermore, as the notation suggests, the differentials determine
extensions in that they reflect the action of the map induced by
$\theta$, namely $\Sq$.  That is, if $x \in H^{n}(\quotient)$ is an
infinite cycle and not a boundary in the spectral sequence, then it
survives to a ``$\Sq$--periodic'' element $\wbar{x}$ in
$H^{n}(\Lambda)$. More precisely, in this situation
$\wbar{\theta}^{i} x$ survives to $E_{\infty}^{n+i,-i}$ for each
$i \geq 0$, and this infinite family corresponds in the abutment to
the family
\[
\{(\Sq)^{i}\wbar{x} \in H^{2^{i} n}(\Lambda) \suchthat i \geq 0 \}.
\]
The presence of a differential $d_{r} \co x \longmapsto y$ from a class
$x \in H^{n,u}(\quotient)$ to a class $y \in H^{n+1,u/2^{r}}(\quotient)$ means
that $r$ ``copies'' of $y$ survive to $E_{\infty}$: the elements
$\wbar{\theta}^{i}y$ survive for $0 \leq i \leq r-1$.  These
elements correspond in the abutment to a finite $\Sq$--family
\[
\{\wbar{y}, \Sq \wbar{y}, (\Sq)^{2} \wbar{y}, \dotsc,
(\Sq)^{r-1} \wbar{y} \},
\]
where $\wbar{y} \in H^{n+1,u/2^{r}}(\Lambda)$.  Thus differentials
coming from the $n$--line in $H^{*}(\quotient)$ give information about
the $\Sq$--action on $H^{n+1}(\Lambda)$.

To compute the $r$th differential on an element $x \in E_{r}$, lift
$x$ to a class in $H^{*}(\quotient)$, represent it by a class in 
$\Lambda$, and take the coboundary.  Since $x$ has survived to
the $E_{r}$--term, this coboundary is in the image of $\theta^{r}$, so
apply $(\theta^{r})^{-1}$.  Project back to a class in
$H^{*}(\quotient)$.  For example, $\lambda_{2} \lambda_{2}
\lambda_{1}$ represents a class at the $E_{1}$--term.  Lift it back to
the class with the same name in $\Lambda$ and take its coboundary.  The
result is $\lambda_{1}^{4}$.  This is in the image of $\theta$, and
applying $\theta^{-1}$ gives $\lambda_{0}^{4}$.  The class
$\lambda_{0}^{4}$ projects to a nonzero class in $H^{*}(\quotient)$.
Thus there is a differential $d_{1}(\lambda_{2} \lambda_{2}
\lambda_{1}) = \lambda_{0}^{4}$.  More generally, there is a
differential
\begin{equation}\label{eqn-d1}
d_{1} (\lambda_{2} \lambda_{2} \lambda_{1}^{n}) = \lambda_{0}^{n+3}
\end{equation}
for $n \geq 1$, reflecting the fact that $\lambda_{0}^{n+3}$
represents a nonzero class in $H^{*}(\Lambda)$, while
$\lambda_{1}^{n+3}=\theta (\lambda_{0}^{n+3})$ is a boundary in
$\Lambda$, and thus $\Sq [\lambda_{0}^{n+3}]=0$ for all $n \geq 1$.

This may not be a good way to compute the cohomology of the Lambda
algebra, but it is a way to study $\Sq$ acting on it.  A little
analysis leads one to believe that most of the cohomology of $\Lambda$
is ``$\Sq$--periodic''.  The word ``most'' in the previous sentence
doesn't have any meaning really, but consider this: given a monomial
$M=\lambda_{n_{1}} \dotsb \lambda_{n_{k}}$ in $\Lambda$ which maps to
a nonzero class in $\quotient$, then at least one $n_{i}$ is even.
Recall that the differential is a derivation, and note that
$d(\lambda_{2m})$ is a sum of terms $\lambda_{i} \lambda_{j}$ with one
of $i$ and $j$ even, the other odd.  Thus applying $d$ to $M$ (in
$\Lambda$) yields a sum of monomials which, before applying the
Adem relations, have the same number of even $\lambda$s as $M$ does.

Furthermore, the Adem relations preserve the parity of the number of
even $\lambda$s.  So there is really only one way to get cohomology
classes in $\quotient$ which are not images of cohomology classes in
$\Lambda$: all terms in the differential must contain an even number
of even $\lambda$s, and the Adem relations on those terms must convert
them into terms with no even $\lambda$s.  The resulting boundary is
zero in the quotient $\quotient$, and hence one gets a cocycle in
$H(\quotient)$ which supports a boundary in the spectral sequence.
Thus such cocycles are precisely those classes in $H(\quotient)$ which
are cocycles but which are not in the image of any cocycle in
$H(\Lambda)$.

Let's look for such things.  First, here is a summary of the previous
paragraphs.

\begin{lemma}\label{lemma-parity}
If a class $x \in H^{*}(\quotient)$ supports a nonzero
differential in the Bockstein spectral sequence~\eqref{eqn-bss}, then
$x$ lifts to a sum of monomials in $\quotient$ which
each contain a positive and even number of even $\lambda$s.
\end{lemma}

Here are two other useful observations in our search.

\begin{lemma}[Wang {\cite[Proposition 1.9]{wang.j*;hopf}}]
\label{lemma-key-fact} 
Given $x \in \quotient$ in positive stem degree, if
$d(x) = 0$ and $x$ is not a boundary, then $x$ is homologous to an
element $y$ with ``odd ending integers'' -- that is, every term in the
admissible expression for $y$ ends in an odd lambda.
\end{lemma}

(Wang's proof works $\quotient$.)

\begin{lemma}[Wang {\cite[Proposition 1.8.3]{wang.j*;hopf}}]
\label{wang-1.8.3}
Let $\lambda_{n_{1}} \dotsb \lambda_{n_{r}} \in \Lambda$ be an
admissible sequence.  Then $(d \lambda_{n_{1}}) \lambda_{n_{2}} \dotsb
\lambda_{n_{r}}$ is a sum of admissible terms with leading integers at
most $n_{1}-1$.
\end{lemma}

\begin{note}\label{note-splitting}
There is an obvious vector space splitting of the short exact sequence
\[
0 \rightarrow \Lambda \xrightarrow{\theta} \Lambda \rightarrow
\quotient \rightarrow 0,
\]
in which one views $\quotient$ as being spanned by the admissible
monomials containing at least one even lambda.  We will often take
this point of view.
\end{note}

The following proposition and theorem are the main results in this
section.

\begin{prop}
Consider the Bockstein spectral sequence~\eqref{eqn-bss}.
\begin{enumerate}
\item There are no nonzero differentials coming from
$H^{1}(\quotient)$, and thus $\Sq$ is injective on $H^{2}\Lambda = 
\ext{A}{2}$.
\item There are no nonzero differentials coming from
$H^{2}(\quotient)$, and thus $\Sq$ is injective on $H^{3}\Lambda =
\ext{A}{2}$.
\end{enumerate}
\end{prop}

Wang proved this result in \cite[Proposition 2.3]{wang.j*;hopf}.

\begin{proof}
(a) \fullref{lemma-parity} implies that any class which supports a
nonzero differential must be in homological degree at least 2.

(b) By \fullref{lemma-parity}, if a class supports a differential,
then it is a sum of terms of the form $\lambda_{2m} \lambda_{2n}$, and
hence is in even stem degree.  Also, \fullref{lemma-key-fact} says
that any cohomology class is cohomologous to a sum of odd-ending
monomials.  Therefore parity considerations imply that any class in an
even stem on the 2--line is cohomologous to a sum of terms of the form
$\lambda_{2i+1} \lambda_{2j+1}$, and hence is zero in $\quotient$.
\end{proof}

\begin{thm}\label{thm-degree-three}
In the Bockstein spectral sequence~\eqref{eqn-bss}, the only
differential emanating from $H^{3}(\quotient)$ is $d_{1} \co \lambda_{2}
\lambda_{2} \lambda_{1} \longmapsto \lambda_{0}^{4}$.  Hence the
kernel of $\Sq$ on $\Ext^{4}$ is spanned by $h_{0}^{4}$.
\end{thm}

\begin{proof}
Suppose that $y \in H^{3}(\quotient)$ supports a differential.  By
\fullref{lemma-key-fact}, we may assume that $y$ is a sum of
admissible monomials each of which ends in an odd lambda, and by
\fullref{lemma-parity}, each term in $y$ has two even lambdas, and
thus $y$ is a sum of terms of the form $\lambda_{\even}
\lambda_{\even} \lambda_{\odd}$.

Following Wang, we write $y$ in the form
\[
y=\lambda_{2n} y_{1} + y',
\]
where $y'$ is a sum of admissible terms with leading term less than
$\lambda_{2n}$ (and hence leading term no larger than
$\lambda_{2n-2}$), and $y_{1}$ is a polynomial in homological degree
2; note that the stem degree of $y_{1}$ is odd.  We also assume that
$y$ is chosen from its cohomology class so that its leading term, in
the lexicographic ordering, is as small as possible.

View $y$ as being an element of $\Lambda$, via the splitting mentioned
in \fullref{note-splitting}.  Since $y$ is a cocycle in $\quotient$, the
boundary of $y$ in $\Lambda$ must consist of all odd terms.  So by
\fullref{wang-1.8.3}, $y_{1}$ must be a cocycle in $\Lambda$.  If
$y_{1}$ were a coboundary, then $y$ would be cohomologous to a class
with smaller leading term, so this can't happen by the minimality
assumption.  Now we appeal to Wang's computation of $H^{2}(\Lambda)$.
By \cite[Proposition 2.4]{wang.j*;hopf}, the only odd ending cocycles
in odd stems are $\lambda_{0} \lambda_{2^{m}-1}$, which in admissible
form is 
\[
\sum_{j=1}^{m-1} \lambda_{2^{m}-2^{j}} \lambda_{2^{j}-1}
= \lambda_{2^{m}-2} \lambda_{1} + \lambda_{2^{m}-4} \lambda_{3} +
\text{other terms}.
\]
Thus $y$ is of the form 
\[
y=(2n, 2^{m}-2, 1) + (2n, 2^{m}-4, 3) + \text{(smaller terms)},
\]
where ``smaller'' is with respect to the lexicographic ordering.
The coboundary of $y$ is 
\[
d (y) = (2n-1, 2^{m}-3, 1, 1) + (2n-1, 2^{m}-4, 2, 1) +
\text{(smaller terms)}.
\]
The first term is all odd, so is zero in $\quotient$, but the second
term cannot be canceled by any of the smaller terms.  Thus it must
not be present, which means that $m$ must be 2: $y$ is of the form 
\[
y=(2n, 2, 1) + \text{(smaller terms)}.
\]
The leading term of each ``smaller term'' is at most $\lambda_{2n-2}$.
We will show that if $n$ is bigger than 1, then the coboundary of $y$
cannot be zero in $\quotient$.  To do this, we will examine the
coboundary of $\lambda_{2n} \lambda_{2} \lambda_{1}$, and find terms
in it which cannot be canceled.  There are three cases, depending on
the congruence class of $n$ mod 4.

\textbf{Case 1} ($n \equiv 0 \mod 4$)\qua
Modulo terms with leading term less than $\lambda_{2n-4}$, here are
the coboundaries of all monomials in the appropriate bidegree:
\begin{align*}
d(\lambda_{2n} \lambda_{2} \lambda_{1}) 
&= 
\lambda_{2n-1} \lambda_{1} \lambda_{1} \lambda_{1} + 
\lambda_{2n-2} \lambda_{1} \lambda_{2} \lambda_{1} + 
\lambda_{2n-4} \lambda_{3} \lambda_{2} \lambda_{1}, \\
d(\lambda_{2n-2} \lambda_{2} \lambda_{3}) 
&= 
\lambda_{2n-2} \lambda_{1} \lambda_{2} \lambda_{1} + 
\lambda_{2n-5} \lambda_{2} \lambda_{2} \lambda_{3}, \\
d(\lambda_{2n-2} \lambda_{4} \lambda_{1}) 
&= 
\lambda_{2n-2} \lambda_{2} \lambda_{1} \lambda_{1} + 
\lambda_{2n-3} \lambda_{3} \lambda_{1} \lambda_{1} + 
\lambda_{2n-3} \lambda_{2} \lambda_{2} \lambda_{1} + 
\lambda_{2n-5} \lambda_{2} \lambda_{4} \lambda_{1}, \\
d(\lambda_{2n-4} \lambda_{6} \lambda_{1}) 
&= 
\lambda_{2n-4} \lambda_{3} \lambda_{2} \lambda_{1} + 
\lambda_{2n-5} \lambda_{5} \lambda_{1} \lambda_{1} + 
\lambda_{2n-5} \lambda_{4} \lambda_{2} \lambda_{1} + 
\lambda_{2n-5} \lambda_{3} \lambda_{3} \lambda_{1} \\
& \quad \quad +  \lambda_{2n-6} \lambda_{4} \lambda_{3} \lambda_{1} + 
\lambda_{2n-6} \lambda_{3} \lambda_{4} \lambda_{1}, \\
d(\lambda_{2n-4} \lambda_{4} \lambda_{3})
&= 
\lambda_{2n-5} \lambda_{2} \lambda_{2} \lambda_{3} + 
\lambda_{2n-6} \lambda_{2} \lambda_{3} \lambda_{3}.
\end{align*}
Any terms which are all odd go to zero in $\quotient$, so we may
ignore them.  Given the monomial $\lambda_{2n}\lambda_{2}
\lambda_{1}$, the only way to cancel $\lambda_{2n-4} \lambda_{3}
\lambda_{2} \lambda_{1}$ in its coboundary is to add $\lambda_{2n-4}
\lambda_{6} \lambda_{1}$, but then there is no way to cancel the term
$\lambda_{2n-5} \lambda_{4} \lambda_{2} \lambda_{1}$.

\textbf{Case 2} ($n \equiv 2 \mod 4$)\qua  If $n=2$, then the leading term of $y$
is $\lambda_{4} \lambda_{2} \lambda_{1}$.  This leads to a permanent
cycle in the spectral sequence: the element
\[
\lambda_{5} \lambda_{1} \lambda_{1} +
\lambda_{4} \lambda_{2} \lambda_{1} + 
\lambda_{2} \lambda_{2} \lambda_{3}
\]
is a cocycle in $\Lambda$, and maps to $\lambda_{4} \lambda_{2}
\lambda_{1} + \lambda_{2} \lambda_{2} \lambda_{3}$ in $\quotient$.

Now assume that $n\geq 6$.  Modulo terms with leading term less than
$\lambda_{2n-5}$, here are the coboundaries of all monomials in the
appropriate bidegree:
\begin{align*}
d(\lambda_{2n} \lambda_{2} \lambda_{1}) 
&= 
\lambda_{2n-1} \lambda_{1} \lambda_{1} \lambda_{1} + 
\lambda_{2n-2} \lambda_{1} \lambda_{2} \lambda_{1} + 
\lambda_{2n-5} \lambda_{4} \lambda_{2} \lambda_{1}, \\
d(\lambda_{2n-2} \lambda_{2} \lambda_{3}) 
&= 
\lambda_{2n-2} \lambda_{1} \lambda_{2} \lambda_{1} + 
\lambda_{2n-5} \lambda_{2} \lambda_{2} \lambda_{3}, \\
d(\lambda_{2n-2} \lambda_{4} \lambda_{1}) 
&= 
\lambda_{2n-2} \lambda_{2} \lambda_{1} \lambda_{1} + 
\lambda_{2n-3} \lambda_{3} \lambda_{1} \lambda_{1} + 
\lambda_{2n-3} \lambda_{2} \lambda_{2} \lambda_{1} + 
\lambda_{2n-5} \lambda_{2} \lambda_{4} \lambda_{1}, \\
d(\lambda_{2n-4} \lambda_{6} \lambda_{1}) 
&= 
\lambda_{2n-4} \lambda_{3} \lambda_{2} \lambda_{1} + 
\lambda_{2n-5} \lambda_{5} \lambda_{1} \lambda_{1} + 
\lambda_{2n-5} \lambda_{4} \lambda_{2} \lambda_{1} + 
\lambda_{2n-5} \lambda_{3} \lambda_{3} \lambda_{1}, \\
d(\lambda_{2n-4} \lambda_{4} \lambda_{3})
&= 
\lambda_{2n-5} \lambda_{2} \lambda_{2} \lambda_{3}.
\end{align*}
The only way to cancel the term $\lambda_{2n-5} \lambda_{4}
\lambda_{2} \lambda_{1}$ in the coboundary of $\lambda_{2n}\lambda_{2}
\lambda_{1}$ is to add $\lambda_{2n-4} \lambda_{6} \lambda_{1}$ to
it, but the coboundary of the resulting sum has a term $\lambda_{2n-4}
\lambda_{3} \lambda_{2} \lambda_{1}$, which cannot be canceled. 

\textbf{Case 3} ($n$ odd)\qua  Assume that $2n > 2$.  Then the
first ``other term'' in the coboundary of $y$ is $\lambda_{2n-3}
\lambda_{2} \lambda_{2} \lambda_{1}$.  To cancel this, the only
appropriate smaller terms are of the form $\lambda_{2n-2} \lambda_{2j}
\lambda_{2k+1}$, where $2j+2k+1 = 5$.  So the terms are
$\lambda_{2n-2} \lambda_{2} \lambda_{3}$ and $\lambda_{2n-2}
\lambda_{4} \lambda_{1}$.  Modulo terms with leading term less than
$\lambda_{2n-3}$, their coboundaries are as follows:
$\lambda_{2n-3}$:
\begin{align*}
d(\lambda_{2n} \lambda_{2} \lambda_{1}) &= \lambda_{2n-1}\lambda_{1}
\lambda_{1} \lambda_{1} + \lambda_{2n-3} \lambda_{2} \lambda_{2}
\lambda_{1}, \\
d(\lambda_{2n-2} \lambda_{2} \lambda_{3}) &= 
\lambda_{2n-2} \lambda_{1} \lambda_{2} \lambda_{1}, \\ 
d(\lambda_{2n-2} \lambda_{4} \lambda_{1}) 
 &= \lambda_{2n-2} \lambda_{2} \lambda_{1} \lambda_{1} +
\lambda_{2n-3} \lambda_{3} \lambda_{1} \lambda_{1} +
\lambda_{2n-3} \lambda_{2} \lambda_{2} \lambda_{1}.
\end{align*}
While we can cancel the term $\lambda_{2n-3} \lambda_{2} \lambda_{2}
\lambda_{1}$, we cannot cancel $\lambda_{2n-2} \lambda_{2}
\lambda_{1} \lambda_{1}$.

So if $y$ supports a differential in the Bockstein spectral sequence,
then $y$ has leading term $\lambda_{2} \lambda_{2} \lambda_{1}$.
There are no smaller terms in the same bidegree, so the potential
cocycle equals its leading term: $y=\lambda_{2} \lambda_{2}
\lambda_{1}$.
\end{proof}

This proof is not ideal.  The technical aspects are cumbersome, and
the proof will not generalize well to higher dimensions.

\subsection[Comments on Sq on the 5--line]{Comments on $\Sq$ on the 5--line}

What about the next degree?  All of the classes in the ideal generated
by $\lambda_{0}^{4}$ must be hit by differentials
\[
d_{1} \co \lambda_{2^{n}-1} \lambda_{2} \lambda_{2} \lambda_{1}
\longmapsto \lambda_{2^{n-1}-1} \lambda_{0}^{4}, \quad \text{for $n \geq
5$}.
\]
Why $n \geq 5$?  Note that the class $\lambda_{2^{n-1}-1} \lambda_{0}^{4}$
is
a coboundary in $\quotient$ when $n=2$, since $\lambda_{1}
\lambda_{0}$ is the coboundary of $\lambda_{2}$.  Thus $\lambda_{1}
\lambda_{0}^{4}$ cannot be the target of a differential in the
spectral sequence.  It is a coboundary when
$n=3$, since $\lambda_{3} \lambda_{0} \lambda_{0}$ is cohomologous to
$\lambda_{1}^{3}$, and $\lambda_{1}^{3} \lambda_{0}$ is the coboundary
of $\lambda_{1} \lambda_{1} \lambda_{2}$.  It is a coboundary when
$n=4$: a calculation shows that it is the coboundary of
\[
\begin{split}
\lambda_{8} \lambda_{0} \lambda_{0} \lambda_{0} + 
\lambda_{6} \lambda_{2} \lambda_{0} \lambda_{0} + 
\lambda_{5} \lambda_{1} \lambda_{2} \lambda_{0} + 
\lambda_{4} \lambda_{4} \lambda_{0} \lambda_{0} + 
\lambda_{4} \lambda_{2} \lambda_{2} \lambda_{0} +
\lambda_{4} \lambda_{1} \lambda_{1} \lambda_{2} \\
+ 
\lambda_{3} \lambda_{3} \lambda_{2} \lambda_{0} + 
\lambda_{2} \lambda_{4} \lambda_{1} \lambda_{1} + 
\lambda_{2} \lambda_{2} \lambda_{4} \lambda_{0} + 
\lambda_{2} \lambda_{2} \lambda_{2} \lambda_{2} + 
\lambda_{1} \lambda_{2} \lambda_{4} \lambda_{1} + 
\lambda_{1} \lambda_{1} \lambda_{2} \lambda_{4}.
\end{split}
\]

As above, suppose that $y=\lambda_{n} y_{1} + y'$ is a cycle in
$\quotient$ with odd-ending integers, in which the leading integers in
the polynomial $y_{1}$ are at most $2n$, and the leading integers in
$y'$ are at most $n-1$.  We may also assume that among classes in
$\quotient$ with odd-ending integers which are homologous to $y$, $y$
has the smallest leading integer.

(We will also view $y$ as an element of $\Lambda$, using the obvious
splitting of the quotient map \fullref{note-splitting}.)

If $n$ is odd, then $y_{1}$ must be a cycle in $\quotient$.  In this
case the boundary of $y$ is
\[
d y = \lambda_{n} (dy_{1}) + (d\lambda_{n}) y_{1} + dy'.
\]
The first term goes to zero in $\quotient$.  We may assume that
$y_{1}$ is nonzero in $\quotient$, and so is not all odd.  We may also
assume that $y_{1}$ is not a boundary.  So either $y_{1}$ is a cycle
in $\Lambda$, or is a non-cycle but has an all-odd boundary.  Suppose
that $y_{1}$ is a cycle in $\Lambda$.  By Wang's results, we know all
of the cycles in $\Lambda^{3}$, and all of the even-dimensional ones
are homologous to an all-odd monomial, and hence to something zero in
$\quotient$.  Therefore we can ignore this case, and we may assume
that $y_{1}$ is a cycle in $\quotient$ which is not the image of a
cycle in $\Lambda$.  By our previous calculation, this means that
$y_{1}=\lambda_{2} \lambda_{2} \lambda_{1}$:
\[
y = \lambda_{2m+1} \lambda_{2} \lambda_{2} \lambda_{1} + y'.
\]  
By an analysis similar to the proof of \fullref{thm-degree-three},
we can show that if $2m+1 \equiv 1 \pmod{4}$, then $y$ cannot be the
leading term of a cocycle.  If $2m+1 \equiv 3 \pmod{4}$, then
$\lambda_{2m+1} \lambda_{2} \lambda_{2} \lambda_{1}$ is the leading
term in the coboundary of $\lambda_{2n+4} \lambda_{2} \lambda_{1}$,
and thus 
\[
\lambda_{2m+1} \lambda_{2} \lambda_{2} \lambda_{1} + y' +
d(\lambda_{2n+4} \lambda_{2} \lambda_{1})
\]
is homologous to $y$ but has smaller leading term.

As a consequence, we may assume that the leading integer $n$ is even.
In this case, then $y_{1}$ must be a cycle in $\Lambda$, and as in the
proof of \fullref{thm-degree-three}, the minimality assumption on
the leading integer of $y$ means that $y_{1}$ cannot be a boundary.  A
result of Wang \cite[Proposition 2.10]{wang.j*;hopf} lists
representatives of all of the cohomology classes of length 3 in an
even stem: they are
\[
\lambda_{2} \lambda_{3} \lambda_{3}, \ \lambda_{0} \lambda_{2^{i}-1}
\lambda_{2^{j}-1}.
\]
Thus we could write down all of the possibilities for the term
$\lambda_{2m} y_{1}$.  This is somewhat complicated, though, because
these are not the admissible forms for the degree 3 cycles, and even
if we had them, left multiplication by $\lambda_{2m}$ could make some
of the resulting terms inadmissible.

Computer calculations in internal degrees up to 82 have found two
phenomena: in this range, there is only one differential other than
those coming from the ideal generated by $\lambda_{2} \lambda_{2}
\lambda_{1}$:
\[
d_{1}(\lambda_{14} \lambda_{6} \lambda_{5} \lambda_{3} +
\text{(smaller terms)}) = \lambda_{6} \lambda_{2} \lambda_{1}
\lambda_{1} \lambda_{1} + \text{(smaller terms)}.
\]
(This reflects the fact that $Ph_{2} \in \ext{A}{5,16}$ is in the
kernel of $\Sq$.)  Also, the source of every differential in this
range except for $\lambda_{2} \lambda_{2} \lambda_{1} \lambda_{1}$ is
cohomologous to a class of the form
\[
\lambda_{2m} \lambda_{6} \lambda_{5} \lambda_{3} + \text{(smaller
terms)}.
\]
This is a relatively small range of dimensions, but one wonders if
these patterns hold for the rest of the 4--line.  Since we don't have
enough evidence to make a good conjecture, we pose questions.

\begin{question}
\begin{enumerate}
\item Is the kernel of $\Sq$ on $\Ext^{5}$ is spanned by the ``obvious''
classes $h_{i}h_{0}^{4}$ with $i \geq 5$, plus the class $Ph_{2}$ in
the 11--stem?
\item Also, consider an element of the form $y=\lambda_{2k} y_{1} + y'
\in \quotient$ which supports a differential in the spectral sequence.
Assume that $y$ is not cohomologous to $\lambda_{2} \lambda_{2}
\lambda_{1} \lambda_{1}$.  Does the lexicographically smallest element
in the cohomology class of $y$ have leading term $\lambda_{2m}
\lambda_{6} \lambda_{5} \lambda_{3}$?
\end{enumerate}
\end{question}

These may be interesting questions, but the methods used in proving
\fullref{thm-main} and \fullref{thm-degree-three} are only going to
get harder in higher filtrations, so some other ideas are needed.

\section{The cocomplete Lambda algebra}\label{sec-cocomplete}

Let $\cocomplete$ be the direct limit of the diagram
\[
\Lambda \xrightarrow{\ \theta\ } \Lambda \xrightarrow{\ \theta\ } \Lambda
\xrightarrow{\ \theta\ } \dotsb.
\]
Call $\cocomplete$ the \emph{cocomplete Lambda algebra}.  For any real
number $r$, write $\zhalf_{>r}$ for the set of elements of $\zhalf$
which are greater than $r$, and similarly for $\zhalf_{\geq r}$.

\begin{prop}
In this proposition, all indices $i$, $j$, $m$, $n$, and $i_{r}$ are
assumed to be in $\zhalf$.
\begin{enumerate}
\item $\cocomplete$ is a $\Z \times \zhalf_{>0}$--graded
$\ztwo$--algebra with generators $\lambda_{n}$ in bidegree $(1,n+1)$
for all $n > -1$.
\item The relations in $\cocomplete$ are generated by
\[
\sum_{i+j=n} \tbinom{i+j}{i} \lambda_{i-1+m} \lambda_{j-1+2m} \quad
\quad \text{for} \ m > 0, \ n \geq 0.
\]
This is the symmetric form of the relations.
\item Alternatively, the relations are generated by 
\[
\lambda_{i} \lambda_{2i+1+n} = \sum_{j \geq 0}
\tbinom{n-j-2^{-N(j,n)}}{j} \lambda_{i+n-j} \lambda_{2i+1+j} 
\quad \quad \text{for} \ i > -1, \ n \geq 0,
\]
where (after \cite[1.7]{llerena-hung}) $N=N(j,n)$ is the least integer
so that $2^{N}j$ and $2^{N}n$ are integers.  This is the admissible
form of the relations.
\item The differential in $\cocomplete$ is given by
\[
d(\lambda_{n-1}) = \sum_{\substack{i+j=n \\ i, j > 0}}
\tbinom{i+j}{j} \lambda_{i-1}\lambda_{j-1} \quad \quad
\text{(symmetric form)}.
\]
\item Alternatively, the differential is given by
\[
d(\lambda_{n}) = \sum_{j > 0}
\tbinom{n-j}{j} \lambda_{n-j}\lambda_{j-1} \quad \quad 
\text{(admissible form)}.
\]
\item The admissible monomials form a basis for $\cocomplete$.  (As in
$\Lambda$, a monomial $\lambda_{i_{1}} \dotsb \lambda_{i_{s}}$ in
$\cocomplete$ is \emph{admissible} if $2i_{r} \geq i_{r+1}$ for $1
\leq r \leq s-1$.)
\item The cohomology of $\cocomplete$ is equal to
$(\Sq)^{-1} \ext{A}{*}$, which in turn is equal to
$\ext{\complete}{*}$, where $\complete$ is the ``complete Steenrod
algebra,'' as studied in \cite{arnon;dickson} and
\cite{llerena-hung}.
\end{enumerate}
\end{prop}

Note that for integers $a$ and $b$, $\tbinom{a}{b} \equiv
\tbinom{2a}{2b} \pmod{2}$, and this allows one to define mod 2 binomial
coefficients for elements of $\zhalf$.  As a consequence, all of the
sums here are finite.  This is essentially because $\cocomplete$ is
constructed as a direct limit; in contrast, the complete Steenrod
algebra $\complete$ is constructed as an inverse limit, and the Adem
relations there are infinite sums -- see \cite[1.7]{llerena-hung}.
Similarly, the admissible monomials in $\complete$ do not span, while
the admissible monomials in $\cocomplete$ do.  The relation with
$\complete$ also explains the terminology ``cocomplete Lambda
algebra.''

\begin{proof}
In general, all of this follows from the colimit description of
$\cocomplete$.  In more detail: for part (b), applying $\theta$ to
the symmetric Adem relation~\eqref{reln-symm} indexed by $m$ and $n$
yields
\begin{align*}
(\text{symmetric Adem relation})_{m,n} &= 
\sum_{i+j=n} \tbinom{i+j}{i} \lambda_{i-1+m} \lambda_{j-1+2m} \\
 & \stackrel{\theta}{\longmapsto} \sum_{i+j=n} \tbinom{i+j}{i}
 \lambda_{2i-1+2m} \lambda_{2j-1+4m} \\
 &= \sum_{2i+2j=2n} \tbinom{2i+2j}{2i} \lambda_{2i-1+2m} \lambda_{2j-1+4m} \\
 &= \sum_{i+j=2n} \tbinom{i+j}{i} \lambda_{i-1+2m} \lambda_{j-1+4m} \\
 &= (\text{symmetric Adem relation})_{2m,2n}.
\end{align*}
In $\Lambda$, one has symmetric Adem relations for all integers $m
\geq 1$ and $n \geq 0$; thus after inverting $\theta$, one needs
relations for all $m, n \in \zhalf$ with $m > 0$ and $n \geq 0$.

Similarly, for part (c), applying $\theta$ to the admissible Adem
relation \eqref{reln-admissible} indexed by $i$ and $n$ yields
\begin{align*}
(\text{admissible Adem } &\text{relation})_{i,n} = 
\lambda_{i} \lambda_{2i+1+n} + \sum_{j \geq 0} \tbinom{n-j-1}{j}
\lambda_{i+n-j} \lambda_{2i+1+j} \\
 & \stackrel{\theta}{\longmapsto} 
\lambda_{2i+1} \lambda_{4i+3+2n} + \sum_{j \geq 0} \tbinom{n-j-1}{j}
\lambda_{2i+1+2n-2j} \lambda_{4i+3+2j} \\
 &= \lambda_{2i+1} \lambda_{4i+3+2n} + \sum_{j \geq 0} \tbinom{2n-2j-2}{2j}
\lambda_{2i+1+2n-2j} \lambda_{4i+3+2j} \\
 &= \lambda_{2i+1} \lambda_{4i+3+2n} + \sum_{j \geq 0} \tbinom{2n-j-2}{j}
\lambda_{2i+1+2n-j} \lambda_{4i+3+j} \\
 &= (\text{admissible Adem relation})_{2i+1,2n}
\end{align*}
In $\Lambda$, one has admissible Adem relations for all non-negative
integers $i$ and $n$, so in $\cocomplete$, one gets admissible Adem
relations for all $i \in \zhalf_{>-1}$ and $n \in \zhalf_{\geq 0}$.
The change of the binomial coefficient from $\tbinom{n-j-1}{j}$ to
$\tbinom{2n-j-2}{j}$ explains the presence of the integer $N(j,n)$ in
the formula.

The two forms of the differentials in $\cocomplete$ are derived
similarly.

Since the admissible monomials form a basis for $\Lambda$, and since
$\theta$ is injective on basis elements, part (f) follows.

Since homology commutes with colimits, part (g) follows.  See also
\cite[5.3]{palmieri;quotient} for the isomorphism between
$(\Sq)^{-1} \ext{A}{*}$ and $\ext{\complete}{*}$.
\end{proof}

The action of $\theta$ on $\cocomplete$ yields an action of an
infinite cyclic group; this action is free in each positive degree in
$\cocomplete$, and trivial in degree zero.  Under this action, the
$\lambda_{n}$'s are partitioned into orbits, and each orbit contains a
unique $\lambda_{n}$ with $n$ an even integer.  Among these,
$\lambda_{0}$ is the only cocycle.  Hence $H^{1}(\cocomplete)$ is
spanned by $\{\theta^{k} (\lambda_{0}) \suchthat k \in \Z \}$.

One would hope that computing $H^{*}(\cocomplete)$ in higher
dimensions would be simpler than computing $H^{*}(\Lambda)$, because
of this symmetry.  We have been unable to take advantage of this so
far, unfortunately.

We remark that $\Lambda$ is a Koszul algebra.  Let $\lie$ be the
``Steenrod algebra for simplicial Lie algebras,'' which is generated
by elements $\SSq^{n}$ and satisfies the usual Adem relations, but has
$\Sq=0$.  Priddy provided a criterion in \cite[5.3]{priddy;koszul}
to check whether an algebra is a Koszul algebra, and he showed in
\cite[8.3--4]{priddy;simplicial-lie} that both $\lie$ and $\Lambda$
satisfy this condition.  Furthermore, he showed in
\cite[9.1]{priddy;simplicial-lie} that $\ext{\lie}{*}$ is isomorphic to
the opposite algebra to $\Lambda$, and pointed out in
\cite[9.4]{priddy;koszul} that 
\[
\ext{\ext{\lie}{*}}{*} \cong \lie.
\]
In other words, $\Lambda$ is the ``Koszul dual'' of $\lie$.

As a consequence, $\cocomplete$ is a Koszul algebra, in a slightly
unconventional sense (since it is $\Z \times \zhalf$--graded, rather
than $\Z \times \Z$--graded).  So one should be able to compute
$\ext{\cocomplete}{*}$ pretty easily; the result should be the
``complete Steenrod algebra for simplicial Lie algebras'' (which is
built from $\lie$ using an inverse limit, just as $\complete$ is built
from $A$).  The details are left for the interested reader.

\appendix

\section[The cohomology of Lambda / theta Lambda through the 14--stem]
{The cohomology of $\Lambda / \theta \Lambda$ through the 14--stem}

\fullref{fig1} contains a table showing the result of hand and computer
calculations of $H^{s,t}(\Lambda / \theta \Lambda)$ up to the 14--stem.

Columns are indexed by stem degree $t-s$, and rows are indexed by
filtration degree $s$.  Each cohomology class is represented by the
leading term of a polynomial representing it, and that leading term
is listed just by giving the subscripts on the lambdas involved;
for example, the entry ``61'' in the 7--stem represents $\lambda_{6}
\lambda_{1} + \lambda_{4} \lambda_{3}$.  The entry ``$0^{n}$'' stands
for $\lambda_{0}^{n}$, while ``$221^{n}$'' stands for $\lambda_{2}
\lambda_{2} \lambda_{1}^{n}$.

This table is complete in the range $t-s \leq 14$ except that it is
missing the higher powers of $\lambda_{0}$ in the 0--stem.  Underlined
classes support differentials in the Bockstein spectral sequence; in
this range, the only differentials are the $d_{1}$'s sending $\lambda_{2}
\lambda_{2} \lambda_{1}^{n}$ to $\lambda_{0}^{n+3}$.

By comparing the well-known computations of $H^{s,t}(\Lambda)$ in this
range with what the Bockstein spectral sequence gives using the visible
differentials, one can conclude that there are no differentials entering
this picture from higher stems.  Therefore the classes $[\lambda_{0}^{i}]$
with $i \geq 4$ are in the kernel of $\Sq$, and for all other classes
$x$ in this range, either $\Sq (x)$ is nonzero or $\Sq (x)$ is out of
this range.

\begin{figure}[ht!]
{\tiny
$$
\begin{turn}{90} 
\text{filtration degree}
\end{turn}
\begin{array}{|c||*{15}{c|}}
\hline
12 & 0^{12} &&&&&&&&&&&&&& \underline{221^{10}} \\ \hline
11 & 0^{11} &&&&&&&&&&&&& \underline{221^{9}}  &\\ \hline
10 & 0^{10} &&&&&&&&&&&& \underline{221^{8}}  && \\ \hline
9  & 0^{9}  &&&&&&&&&&& \underline{221^{7}}  &&& \\ \hline
8  & 0^{8}  &&&&&&&&&& \underline{221^{6}}  &&&& \\ \hline
7  & 0^{7}  &&&&&&&&& \underline{221^{5}}  && 1124111 &&& \\ \hline
6  & 0^{6}  &&&&&&&& \underline{221111}  && 124111 & 224111 &&& 344111 \\ \hline
5  & 0^{5}  &&&&&&& \underline{22111}  && 24111 && 44111 &&& 51233 \\ \hline
4  & 0^{4}  &&&&&& \underline{2211}  & 4111 && 1233 &&&&& 6611 \\ \hline
3  & 0^{3}  &&&&& \underline{221}  && 421 & 431 &&&&&& 671 \\ \hline
2  & 0^{2}  &&& 21 &&&& 61 &&&&&&&  \\ \hline
1  & 0  &&&&&&&&&&&&&&  \\ \hline
0  &&&&&&&&&&&&&&&  \\ \hline\hline
 & 0 & 1 & 2 & 3 & 4 & 5 & 6 & 7 & 8 & 9 & 10 & 11 & 12 & 13 & 14 \\ \hline
\multicolumn{16}{c}{\text{stem degree}}
\end{array}$$}
\caption{Calculations of $H^{s,t}(\quotient)$ up to the 14--stem}
\label{fig1}
\end{figure}

\bibliographystyle{gtart}
\bibliography{link}

\end{document}